\documentclass[12pt,a4paper]{amsart}
\usepackage{amssymb}
\usepackage{amsmath}
\usepackage{amscd}
\usepackage{latexsym}

\headsep=1cm \numberwithin{equation}{section}
\newtheorem{theorem}{Theorem}[section]
\newtheorem{lemma}[theorem]{Lemma}
\newtheorem{proposition}[theorem]{Proposition}
\newtheorem{corollary}[theorem]{Corollary}

\theoremstyle{definition}
\newtheorem{definition}[theorem]{Definition}
\theoremstyle{remark}
\newtheorem{remark}[theorem]{Remark}
\newtheorem{example}[theorem]{Example}

\newtheorem{Definition and Notation}[theorem]{Definition and Notation}
\theoremstyle{Definition and Notation}

\newcommand{\Spec}{\operatorname{Spec}}

\newcommand{\HH}{\operatorname{H}}

\newcommand{\Hom}{\operatorname{Hom}}
\newcommand{\Ext}{\mathrm{Ext}}
\newcommand{\Tor}{\operatorname{Tor}}
\newcommand{\Extindex}{\mathrm{Ext}\text{-}\mathrm{index}}

\newcommand{\Llm}{\lim\limits}

\newcommand{\BN}{\Bbb N}

\newcommand{\fm}{\frak{m}}
\newcommand{\fp}{\frak{p}}

\begin{document}

\title[A note on symmetry in the vanishing of Ext]
{A note on symmetry in the vanishing of Ext}

\author[S. Nasseh]{Saeed Nasseh}
\address{S. Nasseh, Department of Mathematics, Shahid Beheshti University,
Tehran, Iran   -and-   School of mathematics, Institute for research in fundamental sciences (IPM), 19395-5746, Tehran, Iran.}
\email{saeed$\_$naseh@mail.ipm.ir}

\author[M. Tousi]{Massoud Tousi}
\address{M. Tousi, Department of Mathematics, Shahid Beheshti University,
Tehran, Iran  -and-   School of mathematics, Institute for research in fundamental sciences (IPM), 19395-5746, Tehran, Iran.}
\email{mtousi@mail.ipm.ir}

%%%%%%%%%%%%%%%%%%%%%%%%%%%%%%%%%%%%%%
%%%%%%%%%%%%%%%%%%%%%%%%%%%%%%%%%%%%%%
%%%%% abstract %%%%%%%%%%%%%%%%%%%%%%%
%%%%%%%%%%%%%%%%%%%%%%%%%%%%%%%%%%%%%%
%%%%%%%%%%%%%%%%%%%%%%%%%%%%%%%%%%%%%%

\begin{abstract} In \cite{AB} Avramov and Buchweitz proved that for finitely generated modules $M$ and $N$ over a complete intersection local ring $R$, $\Ext^i_R(M,N)=0$ for all $i\gg 0$ implies $\Ext^i_R(N,M)=0$ for all $i\gg 0$. In this note we give some generalizations of this result. Indeed we prove the above mentioned result when (1) $M$ is finitely generated and $N$ is arbitrary, (2) $M$ is arbitrary and $N$ has finite length and (3) $M$ is complete and $N$ is finitely generated.
\end{abstract}

\thanks{This research was in part supported by a grant from IPM (No. 87130214)}

\maketitle

%%%%%%%%%%%%%%%%%%%%%%%%%%%%%%%%%%%%%%
%%%%%%%%%%%%%%%%%%%%%%%%%%%%%%%%%%%%%%
%%%%%% introduction %%%%%%%%%%%%%%%%%%
%%%%%%%%%%%%%%%%%%%%%%%%%%%%%%%%%%%%%%
%%%%%%%%%%%%%%%%%%%%%%%%%%%%%%%%%%%%%%

\section{introduction}
Throughout the paper,  $R$ is assumed to be a commutative Noetherian ring with unity and dim$(R) < \infty$. When $R$ is a local ring, for each $R$-module $M$, $\widehat{M}$ denotes the completion of $M$ with respect to the maximal ideal.

In \cite[Theorem III]{AB} Avramov and Buchweitz proved that for finitely generated modules $M$ and $N$ over a complete intersection local ring $R$, $\Ext^i_R(M,N)=0$ for all $i\gg 0$ implies $\Ext^i_R(N,M)=0$ for all $i\gg 0$. They were interested in determining a class of local rings which satisfy this property. Then Huneke and Jorgensen \cite{HJ} defined a class of Gorenstein local rings, which they called AB rings, and they showed that AB rings satisfy the above mentioned property (see \cite[Theorem 4.1]{HJ}).

Using the notation of \cite{AY}, for given nonzero $R$-modules $M$ and $N$, we define $p^R(M,N)$ to be
$$
p^R (M,N)= \sup \{ i \in \BN  \ | \  \Ext_R^i(M,N) \neq  0 \}.
$$
According to the paper \cite{HJ}, define the Ext-index of the ring $R$,  denoted by $\Extindex (R)$, to be  the supremum of finite values  of  $p^R (M, N)$ for finitely generated  $R$-modules  $M$ and $N$. The authors in \cite{HJ} also called  $R$  an AB ring if it is a Gorenstein local ring of finite Ext-index. Furthermore they showed that the class of AB rings is strictly larger than the class of complete intersection local rings.

In section 2 of this paper we introduce an especial class of AB rings and we show that every complete intersection local ring belongs to this class. Then we show the following theorem:\\
{\bf Theorem A. }{\it Let $R$ be a $d$-dimensional complete intersection local ring. Assume that $M$ and $N$ are two $R$-modules such that $M$ is finitely generated and $N$ is arbitrary. Then
$$
\Ext^i_R(M,N)=0\ \text{for all}\ i\gg 0\Longrightarrow
\Ext^i_R(N,M)=0\ \text{for all}\ i > d.
$$}

In sections 3 we are concerned with the property of symmetry in the vanishing of $\Ext$ over complete intersection local rings when the module appears in the left hand side is not necessarily finitely generated and the right hand side module is finitely generated. As we see in \cite{Si}, it is a general feeling that completeness is a kind of finiteness condition. Therefore in this direction we prove the following theorem:\\
{\bf Theorem B. }{\it Suppose that $R$ is a $d$-dimensional complete intersection local
ring and $M$, $N$ are two $R$-modules. If either $M$ is of finite length and $N$ is arbitrary, or $M$ is finitely generated and $N$ is
complete, then
$$
\Ext^i_R(N,M)=0\ \text{for all}\ i\gg 0\Longrightarrow
\Ext^i_R(M,N)=0\ \text{for all}\ i\gg 0.
$$}

%%%%%%%%%%%%%%%%%%%%%%%%%%%%%%%%%%%%%%
%%%%%%%%%%%%%%%%%%%%%%%%%%%%%%%%%%%%%%
%%%%%%%%% preliminaries %%%%%%%%%%%%%%%%
%%%%%%%%%%%%%%%%%%%%%%%%%%%%%%%%%%%%%%
%%%%%%%%%%%%%%%%%%%%%%%%%%%%%%%%%%%%%%
\section{Preliminaries and Theorem A}

Let $R$ be a Gorenstein local ring and $M$ be a finitely
generated $R$-module. Let $M^*$ denote the dual $R$-module
$\Hom_R(M,R)$. If $M$ is a maximal Cohen-Macaulay (MCM for short) $R$-module, then  there exists a long exact sequence
$$
\mathcal{C}(M): ...\longrightarrow F_2\overset{\partial_2}\longrightarrow
F_1\overset{\partial_1}\longrightarrow
F_0\overset{\partial_0}\longrightarrow
F_{-1}\overset{\partial_{-1}}\longrightarrow
F_{-2}\overset{\partial_{-2}}\longrightarrow
F_{-3}\overset{\partial_{-3}}\longrightarrow ...
$$
of finitely generated free $R$-modules such that $M = Ker {\partial_{-1}}$. Define the
non-negative and negative syzygies of $M$ by $M_i=
Ker\partial_{i-1}$ for every integer $i$. Now we recall \cite[Lemma 1.1]{HJ}, but one should note that the given proof in \cite{HJ} is also true when $N$ is an arbitrary $R$-module, more precisely:

\begin{lemma}\label{first lemma} Let $R$ be a Gorenstein local ring. Suppose that $M$ is an MCM $R$-module and $N$ is an arbitrary $R$-module. Then for fixed $t \geq 3$ and for $1\leq i\leq t-2$ we have
$$
\Ext^i_R(M_{-t},N)\cong \Tor^R_{t-i-1}(M^*,N).
$$
\end{lemma}

\begin{definition} Set $\xi(R)$ to be the supremum of finite values of $p^R (M, N)$ where $M$ and $N$ are $R$-modules and $M$ is finitely generated, i.e.
$$
\begin{aligned}
\xi(R)= \sup \left\{\right. p^R(M,N) \ | \ &\text{$p^R(M,N)<\infty$ where  $M$ is } \\
 &\text{a finitely generated $R$-module} \left.\right\}.
\end{aligned}
$$
We say that the ring $R$ has finite $\xi$ (or is of finite $\xi$) if it satisfies $\xi(R) < \infty$.
\end{definition}

As some obvious properties of this type of rings we point out to the following proposition.

\begin{proposition}\label{first proposition}
\begin{itemize}
\item[]
\item[$(1)$]
Suppose that $(R,\fm)$ is a local ring with $\xi(R) < \infty$. Assume that $x$ is a nonzero divisor on $R$. Then $\xi(R/xR) < \infty$.

\item[$(2)$]
If $R$ is a $d$-dimensional Gorenstein local ring with $\xi(R) < \infty$, then $\xi(R) = d$.

\item[$(3)$]
Every complete intersection local ring is of finite $\xi$.

\item[$(4)$]
Every Gorenstein local ring with finite $\xi$ is an AB ring.

\item[$(5)$]
Suppose that $R$ is a $d$-dimensional Gorenstein local ring with finite $\xi$.
Then for every $\frak p\in \Spec(R)$, $R_{\frak p}$ is of finite
$\xi$.
\end{itemize}
\end{proposition}

\begin{proof}
The proofs of (1), (2) and (3) are completely similar to the proofs of \cite[Propositions 3.3(1), 3.2 and Corollary 3.5]{HJ} respectively. (4) is trivial.

(5) Suppose that $M$ is a finitely
generated $R_{\fp}$-module and $N$ is an arbitrary $R_{\fp}$-module
such that $\Ext_{R_{\fp}}^i(M,N)=0$ for all $i\gg 0$. Write $M =
R_{\fp}y_1+...+R_{\fp}y_t$. Let $M' = Ry_1+...+Ry_t$. We have
$M_{\fp}\cong M\cong M'_{\fp}$. Thus if $\textbf{F.}\longrightarrow
M'\longrightarrow 0$ is a free resolution for $M'$ as an $R$-module,
then $\textbf{F.}\otimes_R{R_{\fp}}\longrightarrow M'_{\fp}\longrightarrow 0$ is a free
resolution for $M'_{\fp}$ as an $R_{\fp}$-module. So, we have
$\Ext^i_R(M',N) \cong \Ext^i_R(M',\Hom_{R_{\fp}}(R_{\fp},N)) \cong \Ext^i_{R_{\fp}}(M,N)$. Therefore
$\Ext^i_R(M',N)=0$ for all $i\gg 0$. By the assumption and (2),
$\Ext^i_R(M',N)=0$ for all $i > d$. Thus
$\Ext^i_{R_{\fp}}(M,N)=0$ for all $i > d$ and this shows that
$\xi(R_{\fp})\leq d$.
\end{proof}

\begin{theorem}\label{main1} Let $R$ be a $d$-dimensional Gorenstein local ring with $\xi(R) < \infty$. Assume that $M$ and $N$ are two $R$-modules such that $M$ is finitely generated and $N$ is arbitrary. Then
$$
\Ext^i_R(M,N)=0\ \text{for all}\ i\gg 0\Longrightarrow
\Ext^i_R(N,M)=0\ \text{for all}\ i > d.
$$
\end{theorem}

\begin{proof}
Let $L$ be the $d$th syzygy of $M$ in a free resolution. We know that $L$ is an MCM $R$-module and $\Ext^i_R(M,N)\cong \Ext^{i-d}_R(L,N)$ for all $i>d$. This shows that $\Ext^i_R(L,N)=0$ for all $i\gg 0$.
Thus for each $t\geq 1$, $\Ext^i_R(L_{-t},N)=0$ for all
$i\gg 0$. Since $\xi(R) < \infty$, then for each $t\geq
1$ and $i > d$, $\Ext^i_R(L_{-t},N)=0$. On the other hand
$\Ext^i_R(L_{-t},N)\cong \Ext^1_R(L_{i-t-1},N)$ for all $i \geq 1$. Thus for each
$t\geq 1$ and $i > d$, $\Ext^1_R(L_{i-t-1},N)=0$. Now by suitable
changing of $i$ and $t$, we will have $\Ext^1_R(L_{-t'},N)=0$ for
each $t'\geq 1$. Therefore by Lemma \ref{first lemma},
$\Tor_{t'-2}^R(L^*,N)=0$ for each $t'\geq 3$.

Therefore if $\textbf{F.}\longrightarrow N\longrightarrow 0$ is
a free resolution for $N$, then
$\textbf{F.}\otimes_R {L^*}\longrightarrow N\otimes_R
{L^*}\longrightarrow 0$ is an exact sequence. Also, $L^*$ is an MCM $R$-module, thus for $i\geq 1$ and every free $R$-module $F$, $\Ext^i_R(F\otimes_R
L^*,R)=0$. So, for $i\geq 1$ we have $\Ext^i_R(N\otimes_R {L^*},R)=\HH^i(\Hom_R(\textbf{F.}\otimes_R {L^*},R))$.
Hence for $i\geq 1$ we get the following isomorphisms
\begin{eqnarray*}
\Ext^i_R(N\otimes_R {L^*},R)&\cong&{\HH^i(\Hom_R(\textbf{F.},\Hom_R(L^*,R)))}\\
&=& {\HH^i(\Hom_R(\textbf{F.},L^{**}))}\\
&\cong& {\HH^i(\Hom_R(\textbf{F.},L))}\\&=&\Ext^i_{R}(N,L).
\end{eqnarray*}
But $R$ is Gorenstein, so $\Ext^i_R(N\otimes_R
{L^*},R)=0$ for $i > d$. Therefore $\Ext^i_{R}(N,L)=0$ for $i > d$. Now since id$(R)=d$ we easily obtain that $\Ext^i_{R}(N,M)=0$ for $i > d$.
\end{proof}

%%%%%%%%%%%%%%%%%%%%%%%%%%%%%%%%%%%%
%%%%%%%%%%%%%%%%%%%%%%%%%%%%%%%%%%%%
%%%%%%%%%%%%%%%%%%%%%%%%%%%%%%%%%%%%
%%%%%%%%%%%%%%%%%%%%%%%%%%%%%%%%%%%%
%%%%%% Proof of theorem B %%%%%%%%%%
%%%%%%%%%%%%%%%%%%%%%%%%%%%%%%%%%%%%
%%%%%%%%%%%%%%%%%%%%%%%%%%%%%%%%%%%%
%%%%%%%%%%%%%%%%%%%%%%%%%%%%%%%%%%%%
%%%%%%%%%%%%%%%%%%%%%%%%%%%%%%%%%%%%

\section{Theorem B}
Let $(R,\fm)$ be a local ring and $E(R/\frak m)$ be the injective
envelope of the residue class field $R/\frak m$. Recall that the
Matlis dual of an $R$-module $T$ is $\Hom_R(T,E(R/\frak m))$ and is
denoted by $T^\vee$. We say that $T$ is Matlis reflexive if
$T^{\vee\vee}\cong T$. Note that if $T$ has finite length, then $T$
is Matlis reflexive. Furthermore we have the following isomorphisms
for $R$-modules $V$ and $W$:
$$
\Tor{_i^R(V,W)}^\vee\cong \Ext^i_R(V,W^\vee)
$$ and
$$
\Ext{^i_R(V,W)}^\vee\cong \Tor_i^R(V,W^\vee) \ \text{when $V$ is
finitely generated.}
$$

\begin{proposition}\label{main2} Suppose that $(R,\frak m)$ is a $d$-dimensional
Gorenstein local ring with finite $\xi$. Then for every $R$-modules $M$ and $N$, where $M$ has finite length and $N$ is arbitrary we have
$$
\Ext^i_{R}(N,M)=0\  \text{for all}\ i\gg 0 \Longrightarrow
\Ext^i_{R}(M,N)=0\  \text{for all}\ i > d.
$$
\end{proposition}
\begin{proof} We have
$$
\Ext^i_{R}(N,M)\cong \Ext^i_{R}(N,M^{\vee\vee})\cong \Tor{_i^R(N,M^\vee)}^{\vee}\cong
\Ext^i_R(M^{\vee},N^{\vee}).
$$
Thus by assumption and Theorem \ref{main1}, $\Ext^i_R(N^{\vee},M^{\vee})=0$ for all $i > d$.
Since $\Ext^i_R(N^{\vee},M^{\vee})\cong
\Tor{_i^R(N^\vee,M)}^{\vee}$, we have $\Tor_i^R(N^\vee,M)=0$ for all
$i > d$. On the other hand $\Tor_i^R(N^\vee,M)\cong \Ext{^i_{R}(M,N)}^\vee$. Therefore $\Ext^i_{R}(M,N)=0$ for all $i > d$.
\end{proof}

By \ref{main1} and \ref{main2} we have the following corollary.

\begin{corollary}\label{Artinian case} Let $R$ be an Artinian Gorenstein local ring with $\xi(R) < \infty$. Assume that $M$ and $N$ are two $R$-modules where $M$ is finitely generated and $N$ is arbitrary. Then
$$
\Ext^i_{R}(N,M)=0\ \text{for all}\ i\gg 0 \Longrightarrow \Ext^i_{R}(M,N)=0\ \text{for all} \ i > 0.
$$
\end{corollary}

\begin{theorem}\label{main3} Suppose that $R$ is a $d$-dimensional Gorenstein local
ring with $\xi(R) < \infty$. Assume that $M$ is a finitely generated $R$-module and $N$ is a
complete $R$-module. Then
$$
\Ext^i_R(N,M)=0\ \text{for all}\ i\gg 0\Longrightarrow
\Ext^i_R(M,N)=0\ \text{for all}\ i\gg 0.
$$
\end{theorem}

To prove this theorem, we need the following definition, remark and lemma.

\begin{definition} Let $(R,\fm)$ be a local ring and $N$ be an arbitrary $R$-module. Let $\tau_N:N\longrightarrow \widehat{N}$ be the natural morphism. We say that $N$ is quasi-complete if $\tau_N$ is surjective and $N$ is separated if $\tau_N$ is injective.
\end{definition}

\begin{remark}\label{recall completion} Suppose that $(R,\fm)$ is a local ring and $N$ is an arbitrary $R$-module. Let $0\longrightarrow K\longrightarrow L\longrightarrow L/K\longrightarrow 0$ be an exact sequence of $R$-modules. From \cite[\S 8]{M}, recall that

(1) $N$ is separated if and only if $\bigcap_n\fm^nN=0$ for all $n\in \Bbb N\cup\{0\}$.

(2) $L/K$ is separated if and only if $K$ is closed in $L$.

(3) Using \cite[1.2, Corollary]{Si} and (2), we get that if $K$ is closed in $L$ and $L$ is quasi-complete then $L/K$ is complete.

Also from \cite[Definition 2.1.11 and Proposition 2.1.12(i)]{St} we have

(4) For every flat $R$-module $F$ there exists a free $R$-submodule $L\subseteq F$ such that the natural injection $\rho: L\longrightarrow F$ is pure (i.e. $\rho\otimes Id_H:L\otimes_RH\longrightarrow F\otimes_RH$ is injective for every $R$-module $H$) and $L$ is dense in the $\fm$-adic topology of $F$ (i.e. $\bigcap_{n\geq 1}(L+\fm^nF)=F$ or $L+\fm^nF=F$ for all $n$). This implies that $L/\fm^nL\cong F/\fm^nF$ for all $n$. Therefore when $F$ is a complete flat $R$-module we have $F\cong \widehat{L}$. In other words every complete flat $R$-module is the completion of a free $R$-module (and conversely, see \cite[2.4]{Si}).
\end{remark}

\begin{lemma}\label{complete} Let $(R,\fm)$ be a local ring and $M$ be a complete $R$-module in $\fm$-adic topology. Suppose that $x\in \frak m$ is a nonzero divisor on both $R$ and $M$. Let
$$
0\longrightarrow T\longrightarrow F\longrightarrow M
\longrightarrow 0
$$
be an exact sequence of $R$-modules where $F$ is a complete flat $R$-module. Then
both $T$ and $T/xT$ are complete in their $\fm$-adic topology.
\end{lemma}

Before proving the lemma, we should remark that $M/xM$ is not necessarily complete, because $xM$ is not necessarily closed in $M$. The following is an  example of A. M. Simon.
\begin{example}\label{example-completion} Let $R = k[[ X,Y,Z ]]$, where $k$ is a field. Put $M_n = R/(XY-Z^n)$ and let $M$ be the completion of $\bigoplus_{n=1}^{\infty}M_n$ as described in \cite[9.4]{Si}. In fact
$$
M=\{(m_n)_{n\geq 1}\in \prod_{n=1}^{\infty}M_n|\ \text{for all}\ s \text{, all but finitely many}\ m_n\ \text{belong to}\ \fm^sM_n\}.
$$
Thus $M\subset \prod_{n=1}^{\infty}M_n$. Note that $X$ is regular on $R$ and $M$. Denote the images of $X$, $Y$, $Z$ in $M_n$ with $x_n$, $y_n$, $z_n$. Let $w_t = (z_1,z_2^2,z_3^3,...,z_t^t,0,...)$ for each $t$. We have that $w_t = X.v_t$, where $(v_t)_i = y_i$ if $i\leq t$ and $(v_t)_i = 0$ otherwise. Thus $w_t\in XM$. The Cauchy sequence $w_t$
has its limit in $M-XM$; indeed we have
$$
\Llm_{t\rightarrow \infty} w_t = (z_1,z_2^2,z_3^3,...,z_t^t,z_{t+1}^{t+1},...)
$$
and $(z_1,z_2^2,z_3^3,...,z_t^t,z_{t+1}^{t+1},...) = X (y_1,y_2,y_3,...,y_t,...)$ which is not in $XM$ because by the above mentioned structure of $M$, $(y_1,y_2,y_3,...,y_t,...)$ is not an element of $M$.
\end{example}

However if $F$ is a complete flat $R$-module, it is mentioned in Remark \ref{recall completion}(4) that $F$ is the completion of a free $R$-module and we observe that $F/xF$ is complete, i.e. $xF$ is closed in $F$. With this in hand we prove Lemma \ref{complete}.

\begin{proof} Since $M$ is complete, $T$ is closed in $F$ and thus
complete (see \cite[1.3, Proposition]{Si}). With our hypothesis, we also have an exact sequence
$$
0\longrightarrow T/xT\longrightarrow F/xF\longrightarrow M/xM\longrightarrow 0.
$$
Thus $xT= T \cap xF$ and $xT$ is closed in $T$ because $T\longrightarrow F$ is continuous.
Consequently by Remark \ref{recall completion}(3), $T/xT$ is complete.
\end{proof}

Now we can give the proof of Theorem \ref{main3}:

\begin{proof} We proceed by induction on $d$. The case $d = 0$ has been
proved in a stronger form in \ref{Artinian case}.

So, assume that $d \geq 1$. Suppose that
$\textbf{P.}\longrightarrow M\longrightarrow 0$ is a free
resolution of $M$. Consider
the short exact sequence $\Lambda: 0\longrightarrow
M_1\longrightarrow P_0 \longrightarrow M\longrightarrow
0$. Since id$P_0 = d$, using the short exact sequence $\Lambda$ and
hypothesis, we have $\Ext^i_{R}(N,M_1)=0$ for all $i\gg 0$.

On the other hand, by \cite[2.5, Proposition]{Si}, there exists a complete flat
resolution $\textbf{F.}\longrightarrow N\longrightarrow 0$ for
$R$-module $N$. By Lemma \ref{complete}, $N_2$ is a complete $R$-module. Also
by \cite[page 85, Corollary 3.2.7]{RG}, for all $j$ we have pd$F_j\leq d$.
Using the exact sequences $\Omega_j: 0\longrightarrow
N_j\longrightarrow F_{j-1}\longrightarrow N_{j-1} \longrightarrow 0$ for $j=1,2$
and the fact that pd$F_{j-1}\leq d$, we obtain
$\Ext^i_{R}(N_2,M_1)=0$ for all $i\gg 0$. Since depth$(R)\geq 1$, there
exists an element $x$ of $\fm$ which is non-zero divisor on $R$, $M_1$
and $N_2$. Thus we have the long exact sequence
$$
\Ext^i_{R}(N_2,M_1)\longrightarrow \Ext^{i+1}_{R}(N_2/xN_2,M_1)
\longrightarrow \Ext^{i+1}_{R}(N_2,M_1)
$$
obtained from the short exact sequence
$$
0\longrightarrow N_2\overset{x.}\longrightarrow N_2\longrightarrow
N_2/xN_2\longrightarrow 0\ \ (\ddagger).
$$

By hypothesis, we have $\Ext^i_{R}(N_2/xN_2,M_1)=0$ for all $i\gg 0$.
Therefore by \cite[Page 140, Lemma 2]{M}, $\Ext^i_{R/xR}(N_2/xN_2,M_1/xM_1)=0$ for all $i\gg 0$.

Now $R/xR$ is a $(d-1)$-dimensional Gorenstein local ring with finite $\xi$ (see Proposition \ref{first proposition}). Also by Lemma \ref{complete}, all $N_i/xN_i$ are complete for $i\geq 2$
and consequently by inductive hypothesis we have
$\Ext^i_{R/xR}(M_1/xM_1,N_2/xN_2)=0$ for all $i \gg 0$. Therefore again by \cite[Page 140, Lemma 2]{M}, $\Ext^i_{R}(M_1,N_2/xN_2)=0$ for all $i \gg 0$. Using again the short exact sequence $(\ddagger)$, we obtain the long
exact sequence
$$
\Ext^i_{R}(M_1,N_2/xN_2)\longrightarrow \Ext^{i+1}_{R}(M_1,N_2)
\overset{x.}\longrightarrow \Ext^{i+1}_{R}(M_1,N_2)\longrightarrow
\Ext^{i+1}_{R}(M_1,N_2/xN_2).
$$

So, we have $\Ext^i_{R}(M_1,N_2) = x \Ext^i_{R}(M_1,N_2)$ for all $i \gg 0$. Therefore by \cite[page 233]{Si}, $\Ext^i_{R}(M_1,N_2)=0$ for all $i \gg 0$. Now, because $R$ is a Gorenstein ring, by \cite[page 79, 3.3.4(ii)]{C}, id$F_j\leq d$. So, using the exact sequences
$\Omega_j$ ($j=1,2$), we have $\Ext^i_{R}(M_1,N)=0$ for all
$i \gg 0$ and using again the exact sequence $\Lambda$, we obtain
$\Ext^i_{R}(M,N)=0$ for all $i \gg 0$.
\end{proof}

As another application of \cite[page 233]{Si} with the same method as above, we close this note by proving the following proposition.

\begin{proposition}\label{main 4} Let $(R,\fm)$ be a $d$-dimensional Gorenstein complete local ring with finite $\xi$. Set
$$
\begin{aligned}
\xi'(R)= \sup \left\{\right. p^R(N,M) \ | \ &\text{$p^R(N,M)<\infty$ where  $M$ is finitely generated} \\
 &\text{and $N$ is arbitrary} \left.\right\}.
\end{aligned}
$$
Then we have $\xi'(R) = d$.
\end{proposition}

\begin{proof} By Corollary \ref{Artinian case} and Theorem \ref{main1}, the claim obviously holds for
$d = 0$. Suppose that $d > 0$. Since id$(R) = d$, there exists an
$R$-module $L$ such that $\Ext^d_{R}(L,R)\neq 0$ and
$\Ext^i_{R}(L,R)=0$ for all $i > d$. Thus $\xi'(R)\geq d$.

Let $M$ be a finitely generated $R$-module and $N$ be an arbitrary
$R$-module such that $\Ext^i_{R}(N,M)=0$ for all $i\gg 0$. Since id$(R)=d$, we can replace $M$ and $N$ by their first syzygies in their $R$-free resolutions. Thus there exists a nonzero divisor $x$ on
$R$, $M$ and $N$. Also using the short exact sequence
$0\longrightarrow M\overset{x.}\longrightarrow M\longrightarrow
M/xM\longrightarrow 0$, we obtain $\Ext^i_R(N,M/xM)=0$ for all $i\gg
0$. Therefore by \cite[page 140, Lemma 2]{M},
$\Ext^i_{R/xR}(N/xN,M/xM)=0\ \text{for all}\ i\gg 0$. Now by
inductive hypothesis we have $\Ext^i_{R/xR}(N/xN,M/xM)=0\ \text{for
all}\ i > d-1$. Therefore, using again the above exact sequence, we
have $\Ext^i_R(N,M) = x \Ext^i_R(N,M)$ for all $i
> d$. But $M$ is a complete $R$-module, so by \cite[page 233]{Si} $\Ext^i_R(N,M)=0$ for all $i
> d$. This shows that $\xi'(R)\leq d$.
\end{proof}

%%%%%%%%%%%%%%%%%%%%%%%%%%%%%%%%%%%
%%%%%%%%%%%%%%%%%%%%%%%%%%%%%%%%%%%
%%%%%%%%%%%%%%%%%%%%%%%%%%%%%%%%%%%
%%%%%%%%%%%%%%%%%%%%%%%%%%%%%%%%%%%
%%%%%%%%%%%%%%%%%%%%%%%%%%%%%%%%%%%
%%%%%%%%%%%%%%%%%%%%%%%%%%%%%%%%%%%
%\newpage
%\vspace{4pt}
\begin{center} \textsc{Acknowledgments}
\end{center}

The authors are grateful to professor A. M. Simon for Example \ref{example-completion} and her other valuable comments on section 3 of this paper.

%%%%%%%%%%%%%%%%%%%%%%%%%%%%%%%%%%%%%%%%%%%%%%%%%%%%%%%%%%%%%%%%%%%%%%%%%%%%%%%%%%%%%%%%%%%%%%%%%%%%%%%%%%%%%%%%%%%%


\begin{thebibliography}{99}

\bibitem{AB}
L.~L.~Avramov and R.~-O.~Buchweitz, {\it Support varieties and cohomology over complete intersections}, Invent. Math. {\bf 142} (2000), 285--318.


\bibitem{AY} T.~Araya and Y.~Yoshino,
{\it Remarks on a depth formula, a grade inequality and a conjecture of Auslander}, Comm. Algebra {\bf 26} (1998), no. 11, 3793--3806.


\bibitem{C}
L.~W.~Christensen, {\it Gorenstein dimensions}, Lecture Notes in Mathematics, {\bf 1747}. Springer-Verlag, Berlin, 2000. viii+204 pp.


\bibitem{RG}
L.~Gruson and M.~Raynaud, {\it Crit$\grave{e}$res de platitude et de projectivit$\acute{e}$. Techniques de "platification" d'un module}, Invent. Math. {\bf 13} (1971), 1--89.


\bibitem{HJ}
C.~Huneke and D.~A.~Jorgensen, {\it Symmetry in the vanishing of Ext over Gorenstein Rings}, Math. Scand. {\bf 93} (2003), 161--184.


\bibitem{M}
H.~Matsumura, {\it Commutative Ring Theory}, Cambridge Studies in Advanced Mathematics, {\bf 8}. Cambridge University Press, Cambridge, 1989. xiv+320 pp.

\bibitem{Si}
A.~M.~Simon, {\it Some homological properties of complete modules}, Math. Proc. Camb. Phil. Soc. {\bf 108} (1990), 231--246.

\bibitem{St}
J.~R.~Strooker, {\it Homological questions in local algebra}, London Math. Soc. Lecture note series, {\bf 145}. Camberidge University Press,  Cambridge, 1990. xiv+308 pp.

\end{thebibliography}
\end{document}